\documentclass{amsart} 

\usepackage{amsmath, amsthm, latexsym, amssymb, amsfonts} 
\usepackage{url}

\urldef{\mailse}\path|als@maildisk.ru| 
\urldef{\mailsf}\path|logo@int-edu.ru| 

\newtheorem{lemma}{Lemma} 
 
\newtheorem*{Theorem A}{Theorem A} 

\newtheorem*{Theorem B}{Theorem B} %-2014dec23-% 

\newtheorem* {seq}{Corollary} 
\newtheorem* {theorem}{Theorem} 
\newtheorem {statement}{Statement} 
\newtheorem {que}{Question} 

\begin{document} 

%\huge 

\title{A Combinatorial Version of the Svenonius Theorem on Definability} 

\author[A.\;L.\;Semenov, \;\,S.\;F.\;Soprunov] 
{A.\;L.\;Semenov, \;\,S.\;F.\;Soprunov
\\\mailse, \mailsf
\\\\Federal Research Center ''Computer Science and Control'' \\of the Russian Academy of Sciences. Moscow (Russia).
}

\begin{abstract} 
The Svenonius theorem describes the (first-order) definability 
in a structure in terms of permutations 
preserving the relations 
of elementary extensions of the structure. 
In the present paper we prove a version of this theorem 
using permutations of sequences over the original structure 
(these are permutations of sequences of tuples 
of the structure elements as well). 
We say that such a permutation $\varphi$ almost preserves a relation 
if for every sequence of its arguments the value of the relation 
on an $n$-th element of the sequence and on its image under $\varphi$ 
coincide for almost all numbers $n.$ 
We prove that a relation is definable in a structure 
iff the relation is almost preserved by all permutations 
almost preserving the relations of the structure. 
This version limits consideration to the original structure only 
and does not refer to any logical notion, such as 
``elementary equivalence''. 


\end{abstract} 

\maketitle 

\bigskip 
{\bf 1. Introduction} \footnote
{This work is supported by the Russian Science Foundation under grant 14-50-00005 with donation of Moscow State University of Education and performed in Steklov Mathematical Institute of Russian Academy of Sciences.}
\bigskip 

Starting in the XIX century the question: 
``How to define something through something?'' 
was considered as a major logical problem. 
From the very beginning the idea of automorphisms 
was associated with the problem. 
The theorem of Svenonius of 1959 plays the role 
% 
%of the Goedel completeness theorem 
% 
of the G\"odel completeness theorem 
for completeness of definability. 
Recent results give more evidence to this view. 
An elaborated survey can be found 
% 
% in [SeSo CSR]. 
in~\cite{se2}. 

The Svenonius theorem (its original version was published in~\cite{sv}) states that un-definability of a relation $R$ in a given structure can be demonstrated always by presenting a permutation of a structure elementarily equivalent to the given structure. The permutation should preserve the relations of the structure and not preserve the relation $R$ and this implies un-definability of the relation. Svenonius (and we as well) uses term ‘permutation’ not ‘automorphism’ to escape ambiguity when considering multiple sets of relations on the common universe. Let us note that the property of “elementary equivalence” in the formulation of  the theorem involves logic. 

We propose to use for the same purposes 
% 
%-2015mar06-% the structure of all sequences over the original structure 
%-2015mar06-% and ``ultimate authomorhisms'' 
%-2015mar06-% (permutations preserving relations everywhere 
%-2015mar06-% but a finite exception). 
permutations on the set of all sequenses of elements 
of the original structure 
and introdue a notion of ``almost preserve''. 
So, we eliminate the need to consider arbitrary structures 
and elementary equivalence and prove our version of the Svenonius theorem. 
%---------------------------- 2015jan15 ----------------------------%
%--------------------------------------------------------------------------

%--------------------------------------------------------------------------
%---------------------------- 2015jan16 ----------------------------%
We start with existing formulations 
and a one-page proof of the classical theorem of Svenonius. 
Then we introduce our main concept of permutation almost 
preserving a relation, prove that it ``fits to'' definability. 
Finally, we prove our version of the Svenonius theorem 
based on this concept.

%A preliminary version of this paper 
%was published in \cite{se1}, 
The result was announced in \cite{se2}. 
%---------------------------- 2015jan16 ----------------------------%
%--------------------------------------------------------------------------

%==========================================================================
%//////////////////////////////////////////////////////////////////////////
%==========================================================================

%--------------------------------------------------------------------------
%---------------------------- 2015jan16 ----------------------------%
%\bigskip 
%{\bf Funding} 
%\nopagebreak 
%\bigskip 

%This work is supported by the Russian Science Foundation under grant 14-50-00005 with donation of Moscow State University of Education %and performed in Steklov Mathematical Institute of Russian Academy of Sciences.
\bigskip 

{\bf Acknowledgments} 
\nopagebreak 
\bigskip 

The authors are thankful to Prof.~Vladimir Uspensky 
for his inspiration in mathematical research, 
to the referees for helpful valuable suggestions, 
and to Yuri Boravlev for technical assistance.

\bigskip 
{\bf 2. Basic definitions} 
\nopagebreak 
\bigskip 
%---------------------------- 2015jan16 ----------------------------%
%--------------------------------------------------------------------------

Let us start with a precise notion of definability. 
% 
%-----2014nov23-----% Let $S$ be some set of relations on a universe $A$ 
Let $S$ be a set of relations on a universe $A$ 
%-----2014nov23-----% and $R$ is a name of a relation on $A$. 
and $R$ be a relation on~$A$. 
%-----2014nov23-----% To define a relation $R$ through $S$ 
%-----2014nov23-----% in a logical language $L$ means: 
%-2015jan22-% To define the relation $R$ in $S$ 
%-2015jan22-% and a logical language $L$ means: 
To {\it define the relation} $R$ {\it in} $S$ 
and a {\it logical language} $L$ means: 

(1) to give names to some relations from $S$ and 

%-----2014dec14-----% (2) to write a formula in the language $L$ using 
%-----2014dec14-----% the given names as extra-logical symbols that 
%-----2014dec14-----% %-----2014nov23-----% is equivalent to $R$ (on $A$). 
%-----2014dec14-----% defines $R$ (on $A$). 

(2) to write a formula that defines $R$ (on $A$) 
in the language $L$ using the given names as extra-logical symbols. 

%-2014dec22-% In this paper $L$ will be the first-order logic 
%-2014dec22-% with equality, we consider countable universe, 
%-2014dec22-% and countable or finite sets of relations. 
% 
%-2014dec31-% In this paper $L$ will be the first-order logic with equality. 
%-2014dec31-% We consider countable universe, 
%-2014dec31-% and countable or finite sets of relations. 
% 
In this paper $L$ will be the first-order logic with equality. 

%--------------------------------------------------------------------------
%---------------------------- 2015jan16 ----------------------------%
%-2015jan22-% We start with existing formulations 
%-2015jan22-% and a one-page proof of the classical theorem of Svenonius. 
%-2015jan22-% Then we introduce our main concept of permutation almost 
%-2015jan22-% preserving a relation, prove that it ``fits to'' definability. 
%-2015jan22-% Finally, we prove our version of Svenonius theorem 
%-2015jan22-% based on this concept. 
%-2015jan22-% 
%-2015jan22-% {\bf 2. Basic definitions} 
%---------------------------- 2015jan16 ----------------------------%
%--------------------------------------------------------------------------

%-2015jan19-% We consider countable universes, 
%-2015jan19-% and countable or finite sets of relations. 

%-2014dec22-% The \emph{(definability) closure} is the 
%-2014dec22-% %-----2014nov23-----% operation of 
%-2014dec22-% extension of a set of relations $S$ 
%-2014dec22-% with all relations definable 
%-2014dec22-% %-----2014nov23-----% through $S$. 
%-2014dec22-% in $S$. 
% 
The \emph{(definability) closure} 
of a set of relations $S$ 
is the extension of it 
with all relations definable in $S$. 
% 
%-2014dec22-% This operation is 
%-2014dec22-% a closure operation in the usual topological 
%-2014dec22-% %-----2014nov23-----% and 
%-2014dec22-% or 
%-2014dec22-% algebraic sense. 
This operation is a closure operation in the usual 
topological or algebraic sense. 
Closed sets of relations we call \emph{definability spaces}, 
the set $S$ is a \emph{base} of the definability closure of $S$. 

%--------------------------------------------------------------------------
%---------------------------- 2015jan19 ----------------------------%
We consider countable universes, 
and countable or finite sets of relations. 
%---------------------------- 2015jan19 ----------------------------%
%--------------------------------------------------------------------------
We call a definability space \emph{countable} %-2014dec22-% 
if it is countable or finite                  %-2014dec22-% 
as a set                                        %-2015mar06-% 
% 
%-2015jan19-% and its domain is countable.                  %-2014dec22-% 
% 
and its universe is countable. %-2014dec22-% %-2015jan19-% 
We assume the Axiom of Choice in our considerations. %-2015jan19-% 

% ´ ª

\medskip 
\textbf{Remark.} 
% 
%-2014dec31-% Our definitions of definability space and other notions are 
%-2014dec31-% %-2014dec22-% "coordinate-free invariant" 
%-2014dec31-% "coordinate-free", "invariant" 
%-2014dec31-% in the sense similar to invariance 
%-2014dec31-% %´coordinate-free invariantª in the sense similar to invariance 
%-2014dec31-% %\ÎÍ coordinate-free invariant\ÔÍ{} 
%-2014dec31-% %in the sense similar to invariance 
%-2014dec31-% of notions for linear spaces, abstract algebras, etc. 
%-2014dec31-% %-2014dec22-% In our case we have properties invariant 
%-2014dec31-% %-2014dec22-% to the choice of names. 
%-2014dec31-% In our case we have properties invariant 
%-2014dec31-% under the change of names. 
% 
Our definitions of definability space and other notions are 
"coordinate-free", "invariant" 
in a sense similar to invariance 
of notions for linear spaces, abstract algebras, etc. 
In our case we have properties invariant 
under the change of names for relations. 
\medskip 

%-2015jan19-% We use following notations: $Dom(f),$ $Im(f)$  %-2014dec22-% 
%-2015jan19-% are the domain and the image of a mapping $f$. %-2014dec22-% 
%-2015jan19-% % 
%-2015jan19-% We denote by $\overline x$ tuples of variables,  %-2014dec22-% 
%-2015jan19-% by $\overline a$ tuples of elements of $A,$ etc. %-2014dec22-% 
% 
We use following notations: $Dom(f),$ $Im(f)$ 
are the domain and the image of a mapping $f,$ 
we denote by $\overline x$ tuples of variables, 
by $\overline a$ tuples of elements of $A,$ etc. 

We denote by $Sym(A)$ the set of all permutations of $A.$ 
% 
%-2014dec31-% A permutation $\varphi$ of $A$ 
%-2014dec31-% \emph{preserves} a relation $R$ iff 
%-2014dec31-% $R(\overline a) \equiv R(\varphi (\overline a)).$ 
% 
A permutation $\varphi$ of $A$ \emph{preserves} a relation $R$ iff 
$R(\overline a) \equiv R(\varphi (\overline a))$ 
for all $\overline a$ from $A.$ 
% 
%-2014dec22-% %-----2014nov23-----% where $\overline a$ is 
%-2014dec22-% for any $\overline a$ being 
%-2014dec22-% a tuple of elements from $A$ 
%-2014dec22-% and $\varphi (\overline a)$ is the tuple of 
%-2014dec22-% %-----2014nov23-----% their images 
%-2014dec22-% images of $\overline a$ 
%-2014dec22-% under permutation $\varphi$. 
% 
%-2015jan19-% A permutation preserves a set of relations $S$ 
%-2015jan19-% if it preserves all relations from $S$, 
%-2015jan19-% a collection of permutations $F$ preserves $S$ 
%-2015jan19-% if every permutation from $F$ preserves~$S.$ 
A permutation {\it preserves a set of relations} $S$ 
if it preserves all relations from $S$, and
{\it a collection of permutations} $F$ {\it preserves} $S$ 
if every permutation from $F$ preserves~$S.$ 

%-2014dec31-% %-----2014nov23-----% With each set of relations $S$ 
%-2014dec31-% With any set of relations $S$ 
%-2014dec31-% we can associate the group $G_S \subseteq Sym(A)$ 
%-2014dec31-% of permutations of the set $A$ preserving $S$. 
%-2014dec31-% %-2014dec22-% It's obvious, that 
% 
%-2015jan19-% With any set of relations $S$ 
%-2015jan19-% we can associate the group $G_S \subseteq Sym(A)$ 
%-2015jan19-% of all permutations of the set~$A$ preserving~$S$. 
With every set of relations $S$ on $A$ 
we can associate the group \mbox{$G_S \subseteq Sym(A)$} 
of all permutations of the set~$A$ preserving~$S$. 
It is obvious, that 
$$S_1 \subseteq S_2 \Rightarrow G_{S_1} \supseteq G_{S_2}$$ 
% 
%-----2014nov23-----% but usually we can't recover a definability space 
but usually we cannot recover a definability space 
from the corresponding subgroup of $Sym(A)$. 

%--------------------------------------------------------------------------
%---------------------------- 2015jan19 ----------------------------%
%-2015jan19-% {2. Ì‡ Ò. 2} 
% 
Up to now we used names for relations implicitly. 
A more standard way is to consider {\it structures:} 
$\langle A, \Sigma, v \rangle,$ 
where $A$ is a universe, 
$\Sigma$ is the set of symbols of an alphabet 
(finite or countable in our case), 
and $v$ is a value (interpretation) 
\mbox{i.\,e.} a function that maps 
% 
%-2015mar06-% any 
each 
symbol from $\Sigma$ 
into a relation over $A.$ 
%---------------------------- 2015jan19 ----------------------------%
%--------------------------------------------------------------------------
% 
%--------------------------------------------------------------------------
%---------------------------- 2015mar06 ----------------------------%
Often the mapping $v$ is omitted. 
In our consideration we do this as well. 
Let \mbox{$M = \langle A, \Sigma, v \rangle$} 
be a structure. 
% 
%-----2014nov23-----% by $tp_\Sigma(\overline b/X)$ 
% 
%-2015mar06-% By $tp_\Sigma(\overline b/X)$ 
%-2015mar06-% we denote the complete type in the signature $\Sigma$ 
%-2015mar06-% %-2015jan20-% of a tuple $\overline b \in M$ over a set $X$ 
%-2015mar06-% of a tuple $\overline b \in A$ over a set $X$ 
%-2015mar06-% %-2014dec24-% i.e. 
%-2015mar06-% \mbox{i.\,e.} 
%-2015mar06-% the set of all formulas 
%-2015mar06-% $\psi(\overline x, \overline a)$ such that 
%-2015mar06-% % 
%-2015mar06-% %-2014dec22-% $M \vDash \psi(\overline b, \overline a), 
%-2015mar06-% %-2014dec22-% \overline a \in X$. 
%-2015mar06-% $$M \vDash \psi(\overline b, \overline a), \quad 
%-2015mar06-% \overline a \in X.$$ 
% 
By $tp_\Sigma(\overline b/X)$ we denote the 
{\it complete type in the signature} 
$\Sigma$ of a 
{\it tuple} $\overline b \in A$ 
{\it over a set} 
$X$ \mbox{i.\,e.} the set of all formulas 
$\psi(\overline x, \overline a)$ such that 
$$M \vDash \psi(\overline b, \overline a), \quad \overline a \in X.$$

%==========================================================================
%//////////////////////////////////////////////////////////////////////////
%==========================================================================

%--------------------------------------------------------------------------
%---------------------------- 2015jan20 ----------------------------%
\bigskip 
{\bf 3. The Svenonius theorem} 
\nopagebreak 
\bigskip 
%---------------------------- 2015jan20 ----------------------------%
%--------------------------------------------------------------------------

%-----2014nov23-----%  Î‡ÒÒË˜ÂÒÍ‡ˇ ÚÂÓÂÏ‡ Svenonius ÒÙÓÏÛÎËÓ‚‡Ì‡ ‚ 
%-----2014nov23-----% \cite{sv} ÒÎÂ‰Û˛˘ËÏ Ó·‡ÁÓÏ: 
Svenonius formulated in \cite{sv} his classical theorem as follows: 

\begin{Theorem A} 
Let $S$ be any elementary system with the predicate constants 
%-2014dec23-% $R, Q, T_1, T_2,\dots,$ 
$R,$ $Q,$ $T_1,$ $T_2,$ $\dots,$ 
such that no disjunction of explicit 
definitions of $Q$ in terms of $R$ is provable in $S$. 
Then there is a model $M$ of $S$ with a permutation $\Phi$ 
which preserves $R$ but not $Q$. 
\end{Theorem A} 

The Svenonius theorem can be formulated in a ``less syntactical'' way 
than the original one 
%-2015mar09-% (compare with Corollary 10.5.2 in \cite{Hodges}). 
(see the Corollary 10.\;\!5.\;\!2 in \cite{Hodges}). 

%-2015jan20-% Theorem B. 
\begin{Theorem B} 
Let $\Sigma$ be a signature, $R$ a symbol, 
$\langle A, \Sigma\cup\{R\}, v \rangle$ a~structure. 
% 
%-2015mar09-% The relation $v(R)$ is not definable in 
%-2015mar09-% $\langle A, \Sigma, v \rangle$ if: 
%
If the relation $v(R)$ is not definable in 
$\langle A, \Sigma, v \rangle$ then 
for every structure $\langle A, \Sigma', v'\rangle,$ 
%\quad 
\; 
$\Sigma' \supseteq \Sigma\cup\{R\}$ 
there exists a structure 
$\langle B, \Sigma', w\rangle$ 
elementary equivalent to 
$\langle A, \Sigma', v'\rangle$ 
and a permutation of $B$ that preserves $w(\Sigma)$ 
and does not preserve $w(R).$ 
\end{Theorem B} 

%-2015mar09-% \medskip 
\textbf{Remark.} 
%-2015mar09-% Note, that both in the Theorem A and Theorem B formulation 
Note, that both in the formulations of Theorem A and Theorem B 
three spaces occur: a definability space generated by $\Sigma,$ 
its extension with a relation and further extension with more relations. 
It is not a superfluous generalization of a more simple formulation 
of the Corollary (see below). 
In fact, we can have a large space generated by $\Sigma'$ 
and its subspace generated by $\Sigma,$ 
and want to know does $R$ extend $\Sigma$ 
% 
%-2015mar09-% it 
% 
properly. 
% 
%-2015mar09-% This is situation of 
In the situation we have 
so called reducts, for example, 
reducts of $\langle {\mathbb Q}, > \rangle$ 
% 
%-2015mar09-% (see [ ]). 
(see \cite{mac}). 
Here $\Sigma' = \{>\}.$ 
Proofs in this situation use the structure of  
$\langle {\mathbb Q}, > \rangle$ 
% 
%-2015mar09-% 
and $>$ plays the role of $\Sigma'.$ 
\medskip 

%-2015jan20-% Proof. 
%-2015mar09-% \begin{proof} 
\noindent{\bf Proof of Theorem B.} 
Let the relation $v(R)$ be not definable in 
$\langle A, \Sigma, v\rangle.$ 
%---------------------------- 2015jan20 ----------------------------%
%--------------------------------------------------------------------------
% 
% 
We shall construct 
(by using back-and-forth argument) %-2014dec31-% 
% 
%-2015jan22-% countable chains 
%-2015mar09-% countable chain 
a countable chain 
$$M \preccurlyeq M_0 \preccurlyeq \dots 
\preccurlyeq M_i \preccurlyeq \dots$$ 
of 
%-2015jan22-% $M_i = \langle A_i,\Sigma \rangle$ 
%-2015mar09-% $M_i = \langle A_i,\Sigma', v_i \rangle$ 
$M_i = \langle A_i,\Sigma', v_i \rangle$ -- 
countable elementary extensions of $M$ and 
$$ \varphi_0 \subseteq \varphi_1 \subseteq \dots 
\subseteq \varphi_i \subseteq \dots $$ 
of finite partial mappings 
$\varphi_i \colon A_i \to A_i.$ 
% 
%-2015jan22-% Then we take as 
%-2015jan22-% ${S^+}'$ 
%-2015jan22-% the set of relations definable in 
%-2015jan22-% $\bigcup_i M_i,$ 
%-2015jan22-% and 
%-2015jan22-% prove that 
%-2015jan22-% $\bigcup_i \varphi_i$ 
%-2015jan22-% is the permutation on the universe of $\bigcup_i M_i$ 
%-2015jan22-% preserving $S'$ and not preserving $R'$. 
% 
Then we shall take 
%--------------------------------------------------------------------------
%---------------------------- 2015jan20 ----------------------------%
%-2015jan22-% {4 Ì‡ Ò. 4} 
$B = \bigcup_i A_i,$ 
\, %-2015mar09-% 
$w = \bigcup_i v_i,$ 
\, %-2015mar09-% 
$\varphi = \bigcup_i \varphi_i$ 
\, %-2015mar09-% 
% 
%-2015mar09-% and prove that $\varphi$ is the required permutation. 
and prove that $\varphi$ is a required permutation. 
% 
%---------------------------- 2015jan20 ----------------------------%
%--------------------------------------------------------------------------
% 
% 
%--------------------------------------------------------------------------
%---------------------------- 2015jan20 ----------------------------%
We enumerate the countable set $\bigcup_i A_i$ beforehand. 
%---------------------------- 2015jan20 ----------------------------%
%--------------------------------------------------------------------------

%--------------------------------------------------------------------------
%---------------------------- 2015mar09 ----------------------------%
Let us denote \mbox{$\Sigma \cup \{R\}$} by $\Sigma^*.$ 
%---------------------------- 2015mar09 ----------------------------%
%--------------------------------------------------------------------------

First of all, let us 
%-2014dec31-% say that $n$ is the number of arguments of $R$ and 
consider the type 
$$ 
\{R(\overline x) \not \equiv R(\overline y)\} \; \cup \; 
\left\{Q(\overline x) \equiv Q(\overline y) 
\left| %\mid 
\text{ \parbox{3cm}{\normalsize \small % 
for all formulas $Q$ in the signature $\Sigma$ 
with the same number of arguments as $R$ 
}}\right.\right\}. 
$$ 

By 
% 
%-2014dec31-% the 
% 
the compactness theorem this type is consistent, 
otherwise the relation $R$ would be definable 
through a finite collection of 
% 
%-2014dec31-% $Q_1,$ $\dots,$ $Q_n$. 
%-2015jan22-% $Q_i.$ 
$Q.$ 
So there is an elementary extension $M_0$ 
of the structure $M$ and a tuple 
\mbox{$\overline a, \overline b \in M_0$} 
which realizes this type. 
% 
%-2015jan22-% We set 
%-2015jan22-% \mbox{$\varphi_0(\overline a) = \overline b$.} 
% 
%-2015mar09-% We set 
%-2015mar09-% \mbox{$\varphi_0(\overline a) = \overline b$,} 
%-2015mar09-% \, 
%-2015mar09-% \mbox{$\Sigma^* = \Sigma \cup \{R\}.$} 
We set 
\mbox{$\varphi_0(\overline a) = \overline b$.} 
% 
%-2014dec31-% Note that $\varphi_0$ does not preserve $R$, 
%-2014dec31-% so $\bigcup_i \varphi_i$ does not preserve $R$ as well. 

%-2014dec31-% Then we use back-and-forth argument 
%-2014dec31-% to construct the needed sequence 
%-2014dec31-% $$ 
%-2014dec31-% M \preccurlyeq M_0 \preccurlyeq \dots 
%-2014dec31-% \preccurlyeq M_i \preccurlyeq \dots. 
%-2014dec31-% $$ 
%-2014dec31-% On every step we produce 
%-2014dec31-% $M_{i+1}$ and $\varphi_{i+1}$ from $M_i,$~$\varphi_i$. 

%-2015jan22-% We have 
%-2015jan22-% % $$ 
%-2015jan22-% % tp_\Sigma(Dom(\varphi_0)/\varnothing) = 
%-2015jan22-% % tp_\Sigma(Im(\varphi_0)/\varnothing) 
%-2015jan22-% % $$ 
%-2015jan22-% $$ 
%-2015jan22-% tp_\Sigma\bigl(Dom(\varphi_0)/\varnothing\bigr) = 
%-2015jan22-% tp_\Sigma\bigl(Im(\varphi_0)/\varnothing\bigr) 
%-2015jan22-% $$ 
%-2015jan22-% % 
%-2015jan22-% %-2014dec31-% and prove that 
%-2015jan22-% and will prove that 
%-2015jan22-% % 
%-2015jan22-% % $$ 
%-2015jan22-% % tp_\Sigma(Dom(\varphi_i)/\varnothing) = 
%-2015jan22-% % tp_\Sigma(Im(\varphi_i)/\varnothing) 
%-2015jan22-% % $$ 
%-2015jan22-% $$ 
%-2015jan22-% tp_\Sigma\bigl(Dom(\varphi_i)/\varnothing\bigr) = 
%-2015jan22-% tp_\Sigma\bigl(Im(\varphi_i)/\varnothing\bigr) 
%-2015jan22-% $$ 
%-2015jan22-% by induction. 
% 
We have 
$$ 
tp_{\,\Sigma^*}\!\bigl(Dom(\varphi_0)/\varnothing\bigr) = 
tp_{\,\Sigma^*}\!\bigl(Im(\varphi_0)/\varnothing\bigr) 
$$ 
and will prove that 
$$ 
tp_{\,\Sigma^*}\!\bigl(Dom(\varphi_i)/\varnothing\bigr) = 
tp_{\,\Sigma^*}\!\bigl(Im(\varphi_i)/\varnothing\bigr) 
$$ 
by induction. 
% 
% 
%-2014dec31-% ------------------ 
Note that $\varphi_0$ does not preserve $R$, 
so $\bigcup_i \varphi_i$ does not preserve $R$ as well. 
%-2014dec31-% ------------------ 

%-2015jan22-% On even steps 
On even steps $i \geqslant 0$ 
we take the first item $a$ in $A_i$ 
that is not included in the domain 
of $\varphi_i.$ 
Then we choose an elementary extension 
$M_{i+1}$ of $M_i$ such that 
% $$ 
% tp_\Sigma(a/Dom(\varphi_i)) = tp_\Sigma(b/Im(\varphi_i)) 
% $$ 
% 
%-2015jan22-% $$ 
%-2015jan22-% tp_\Sigma\bigl(a/Dom(\varphi_i)\bigr) = 
%-2015jan22-% tp_\Sigma\bigl(b/Im(\varphi_i)\bigr) 
%-2015jan22-% $$ 
$$ 
tp_{\,\Sigma^*}\!\bigl(a/Dom(\varphi_i)\bigr) = 
tp_{\,\Sigma^*}\!\bigl(b/Im(\varphi_i)\bigr) 
$$ 
for some $b \in A_{i+1}$. 
(Note that by definition of type, 
equalities belong to the types and %-2015jan22-% 
% 
%-2015mar09-% $b$ is different 
so $b$ is different 
from elements in the image of the partial 
% 
%-2015jan22-% automorphism.) 
mapping.) 
Finally, we set 
$$ 
\varphi_{i+1}=\varphi_i \cup \{\langle a,b \rangle\}. 
$$ 

%-2015mar09-% On odd steps we take the first item $b \in A_i$ 
On odd steps $i$ we take the first item $b \in A_i$ 
that is not in the image of $\varphi_i.$ 
Then we choose an elementary extension 
$M_{i+1}$ of $M_i$ such that 
% 
%-2015jan22-% $$ 
%-2015jan22-% tp_\Sigma\bigl(a/Dom(\varphi_i)\bigr) = 
%-2015jan22-% tp_\Sigma\bigl(b/Im(\varphi_i)\bigr) 
%-2015jan22-% $$ 
$$ 
tp_{\,\Sigma^*}\!\bigl(a/Dom(\varphi_i)\bigr) = 
tp_{\,\Sigma^*}\!\bigl(b/Im(\varphi_i)\bigr) 
$$ 
for some 
% $a \in M_{i+1}$. 
\mbox{$a \in M_{i+1}$.} 
(Note that $a$ is not included in the domain 
of the current partial automorphism.) 
Now we set 
$$ 
% \varphi_{i+1}=\varphi_i \cup \{<a,b>\}. 
\varphi_{i+1}=\varphi_i \cup \{\langle a,b \rangle\}. 
$$ 

%-2015jan22-% %-2014dec31-% Note that 
%-2015jan22-% So, 
%-2015jan22-% $\varphi_{i+1}$ preserves 
%-2015jan22-% all relations from $S$ and 
%-2015jan22-% $$ 
%-2015jan22-% tp_\Sigma\bigl(Dom(\varphi_{i+1})/\varnothing\bigr) = 
%-2015jan22-% tp_\Sigma\bigl(Im(\varphi_{i+1})/\varnothing\bigr). 
%-2015jan22-% $$ 
%-2015jan22-% %-2014dec31-% $\bigcup_i \varphi_i$ is the permutation 
%-2015jan22-% %-2014dec31-% of the universe. 
% 
So, for all $i$ the partial mapping 
$\varphi_{i+1}$ preserves all relations from $S$ and 
$$ 
tp_{\,\Sigma^*}\!\bigl(Dom(\varphi_{i+1})/\varnothing\bigr) = 
tp_{\,\Sigma^*}\!\bigl(Im(\varphi_{i+1})/\varnothing\bigr). 
$$ 

%-2015jan22-% Now it is easy to see that $\bigcup_i M_i$, 
%-2015jan22-% $\bigcup_i \varphi_i$ 
%-2015jan22-% are needed elementary extension and permutation. 
%-2015jan22-% \end{proof} 
% 
Now it is easy to see that $\bigcup_i M_i$, 
\, 
$\bigcup_i \varphi_i$ 
are the needed elementary extension and (total) permutation. 
% 
%-2015mar09-% \end{proof} 
%\hfill $\Box\mathsurround=-4pt$ 
%\hfill $\Box$\hspace{-\mathsurround}% 
\hfill $\Box\mathsurround=0pt$% 
%\medskip 
\bigskip 

In the case $\Sigma' = \Sigma \cup \{R\}$ 
we have a simplified version of the theorem. 

%-2015jan22-% Corollary. 
\begin{seq} 
Let $\Sigma$ be a signature, $R$ a symbol, 
$\langle A, \Sigma \cup \{R\}, v \rangle$ a~structure. 
If the relation $v(R)$ is not definable in 
$\langle A, \Sigma, v \rangle$ 
then there exists\- a structure 
$\langle B, \Sigma \cup \{R\}, w \rangle$ 
elementary equivalent to it 
and a~permutation of $B$ that preserves $w(\Sigma)$ 
and does not preserve $w(R).$ 
\end{seq} 
%---------------------------- 2015jan22 ----------------------------%
%--------------------------------------------------------------------------

The Svenonius theorem 
% 
%-2015mar09-% in its last formulation above %-2014dec31-% 
in its last formulation above (Theorem B) 
reduces the question of the first-order 
(un-)definability to a question about automorphisms 
of structures elementarily equivalent to the original one. 
So, indirectly it refers to first-order logic again. 
% 
%-2014dec24-% In this paper we find a way to eliminate logic and 
%-2014dec24-% produce a purely combinatorial form of Svenonius theorem. 
%
The notion of elementary equivalence can be characterized 
%-2015mar09-% in terms of Ehrenfeucht -- Fra\"{i}ss\'{e} games 
in terms of Ehrenfeucht~\mbox{--~Fra\"{i}ss\'{e}} games 
%-2014dec24-% (e.g. \cite{Hodges}) 
(\mbox{e.\,g.} \cite{Hodges}) 
or according to a well-known result of Keisler and Shelah 
in terms of ultrapowers 
%-2014dec24-% (e.g. \cite{ult}), 
%-2014dec24-% (e.g. \cite{ult}). 
(\mbox{e.\,g.} \cite{ult}). 
%-2014dec24-% so the Svenonius theorem 
So the Svenonius theorem 
allows us to describe the first-order definability 
% 
%-2015mar09-% in combinatorial terms of games 
in ``combinatorial'' terms of games 
or ultrapowers and permutations. 
%-2014dec24-% But we propose a more direct and simple construction. 
We propose a more direct and simple construction 
to eliminate logic and 
% 
%-2015jan22-% produce a purely combinatorial form of Svenonius theorem. 
produce a combinatorial form of the Svenonius theorem.

%==========================================================================
%//////////////////////////////////////////////////////////////////////////
%==========================================================================

%--------------------------------------------------------------------------
%---------------------------- 2015jan23 ----------------------------%
\bigskip 
{\bf 4. The main result} 
\nopagebreak 
\bigskip 

Let $A$ be a countable universe, $\mathbb{N}$ 
-- the set of all natural numbers. 
By $\mathcal F$ we denote the set of everywhere defined functions 
$f \colon \mathbb{N} \to A$ 
\mbox{i.\,e.} $A^\mathbb{N}$ 
% 
%-2015mar09-% (the set of all infinite sequences of $A$). 
(the set of all infinite sequences over $A$). 
If $P$ is an $n$-ary relation on $A$ 
and $\varphi$ 
is a (partial) mapping $\mathcal F \to \mathcal F$ 
%is a permutation of $\mathcal F$ 
% 
then we say that $\varphi$ \emph{almost preserves} $P$ if 
$$ 
\left\{i \; \mid \; 
P\bigl(f_1(i),\dots,f_n(i)\bigr) \not \equiv 
P\bigl(\varphi(f_1)(i),\dots,\varphi (f_n)(i)\bigr) 
\right\} 
$$ 
is finite for each $f_1,$ $\dots,$ $f_n$ in $Dom(\varphi)$. 
We say that $\varphi$ \emph{almost preserves} 
a definability space $S$ if it almost preserves every $P$ in $S$. 
% 
%Such $\varphi$ can be called ``ultimate automorphisms''. 
A permutation $\varphi$ which almost preserves 
a definability space can be called ``ultimate automorphisms''.
%--------------------------------------------------------------------------
%---------------------------- 2015jan23 ----------------------------%
The following statement justifies our concept 
of almost preservation in the context of definability. 
%---------------------------- 2015jan23 ----------------------------%
%--------------------------------------------------------------------------

\begin{statement} 
%-2014dec24-% ≈ÒÎË $\Sigma$ -- Ò˜ÂÚÌÓÂ 
%-2014dec24-% ÏÌÓÊÂÒÚ‚Ó ÓÚÌÓ¯ÂÌËÈ Ì‡ $A$ Ë ÓÚÌÓ¯ÂÌËÂ $R$ Ì‡ $A$ 
%-2014dec24-% ÓÔÂ‰ÂÎËÏÓ ‚ $\Sigma$, ÚÓ Î˛·‡ˇ ÔÂÂÒÚ‡ÌÓ‚Í‡, 
%-2014dec24-% ÔÓ˜ÚË ÒÓı‡Ìˇ˛˘‡ˇ ‚ÒÂ ÓÚÌÓ¯ÂÌËˇ ËÁ $\Sigma$ ÔÓ˜ÚË 
%-2014dec24-% ÒÓıˇÌˇÂÚ ÓÚÌÓ¯ÂÌËÂ $R$. 
% 
%-2014dec31-% Let $S$ be a countable set of relations, 
Let $S$ be a countable set of relations on $A$, and suppose that  
\/ %-2015jan23-% 
$R$ is definable in $S.$ 
% 
%-2015jan23-% Then any permutation that almost 
%-2015mar09-% Then every permutation that almost 
Then every permutation of $A^\mathbb{N}$ that almost 
% 
%-2014dec31-% preserve $S$ almost preserve $R.$ 
preserves $S$ almost preserves~$R.$ 
\end{statement}

\begin{proof}
Let $\varphi$ be a permutation on $\mathcal{F}$, 
and suppose that $\varphi$ almost preserves all relations from 
% 
%-2014dec24-% $\Sigma$, 
$S$, 
and for some 
% 
%-2014dec24-% $f_1,\dots,f_n$ 
$f_1,$ $\dots,$ $f_n$ 
in $\mathcal{F}$ the set 
%-2014dec24-% $$ 
%-2014dec24-% V=\{i | 
%-2014dec24-% R(f_1(i),\dots,f_n(i)) \not \equiv 
%-2014dec24-% R(\varphi(f_1)(i),\dots,\varphi (f_n)(i)) \} 
%-2014dec24-% $$ 
$$ 
V= 
\left\{i \; \mid \; 
R\bigl(f_1(i),\dots,f_n(i)\bigr) \not \equiv 
R\bigl(\varphi(f_1)(i),\dots,\varphi (f_n)(i)\bigr) 
\right\} 
$$ 
is infinite. 
% 
%-2014dec31-% Consider 

Let $U$ be a nonprincipal ultrafilter on $\mathbb{N}$ 
% 
%-2014dec24-% $V \in U$. 
%-2015mar09-% \; \mbox{$V \! \in U$}. 
with \mbox{$V \! \in U$}.
Let $M$ be 
% 
%-2015mar09-% $\langle A,S \rangle,$ 
the structure  with universe $A,$ set of names $S$ 
and the given values for $S$. 
Let  $H$ be the
% 
%-2015mar09-% an ultrapower $\prod_U M$ where 
ultrapower \mbox{$\prod_U M$} where 
$U$ is a nonprincipal ultrafilter on $\mathbb{N}$, 
% 
%-2014dec24-% $V \in U$. 
%-2015mar09-% \; \mbox{$V \! \in U$}. 
\,\mbox{$V \! \in U$}. 
% 
% 
% 
% 
% 
% -2014dec25- % 
For $a \in \mathcal{F}$ we denote by $\widetilde{a}$ 
the equivalence class of 
% 
%-2014dec25-% $\{f | \{i | f(i) = a(i)\} \in U\}$ 
$$ 
\bigl\{ 
f \; \mid \; \{i \mid f(i) = a(i)\} \in U 
\bigr\} 
$$ 
% 
%-2015mar09-% in $\prod_U M.$ 
in $H.$ 
We define the mapping 
% 
%-2015mar09-% $$\varphi_U \colon \prod_U M \to \prod_U M$$ 
$\varphi_U \colon H \to H$ 
% 
% 
%-2014dec24-% such that 
%-2014dec24-% $\varphi_U(\widetilde{a})=\widetilde{b}$ 
%-2014dec24-% if 
%-2014dec24-% $\varphi(a) \in \widetilde{b}$ 
%-2014dec24-% for some 
%-2014dec24-% $a \in \widetilde{a}$. 
% 
% 
such that 
\, \mbox{$\varphi_U (\widetilde{a}) \mathop{=} \widetilde{b}$} \, 
if 
\, \mbox{\mbox{$\varphi(a) \mathop{\in} \widetilde{b}$}} \, 
for some 
\, \mbox{$a \mathop{\in} \widetilde{a}$}. 
% 
% 
%-2015jan01-% Note, that for each 
Let us note, that because we consider structures with equality, for each 
% 
%-2014dec24-% $f_1, f_2 \in \mathcal{F}$ 
$f_1, f_2 \mathop{\in} \mathcal{F}$ 
% 
%-2015jan01-% holds 
the set 
% 
% 
%-2014dec24-% $$\{i | (f_1(i) = f_2(i)) \equiv 
%-2014dec24-% (\varphi(f_1)(i) = \varphi (f_2)(i))\}$$ 
% 
$$ 
\left\{ 
i \; \; \left| \; \; %\mid %| 
\bigl(f_1(i) = f_2(i)\bigr) \not \equiv 
\bigl(\varphi(f_1)(i) = \varphi (f_2)(i)\bigr) 
\right. 
\right\} 
$$ 
is finite. 
So 
% 
% 
%-2015mar09-% $$ 
%-2015mar09-% \prod_U M \vDash \widetilde{f_1} = \widetilde{f_2} 
%-2015mar09-% %-2015jan01-% $$ 
%-2015mar09-% % 
%-2015mar09-% %-2015jan01-% iff 
%-2015mar09-% \;\; \Leftrightarrow \;\; %-2015jan01-% 
%-2015mar09-% % 
%-2015mar09-% %-2015jan01-% $$ 
%-2015mar09-% \prod_U M \vDash \widetilde{\varphi}(f_1) = \widetilde{\varphi} 
%-2015mar09-% (f_2) 
%-2015mar09-% $$ 
%-2015mar09-% % 
%-2015mar09-% and the mapping $\varphi_U$ is 
%-2015mar09-% a %-2015jan01-% 
%-2015mar09-% well-defined permutation of 
%-2015mar09-% % 
%-2015mar09-% % 
%-2015mar09-% % 
%-2015mar09-% %-2014dec25-% $\prod_U M$. 
%-2015mar09-% $\prod_U M$ 
%-2015mar09-% (by $\widetilde{\varphi}(a)$ 
%-2015mar09-% we denote the equivalence class 
%-2015mar09-% of~$\varphi(a)$). 
% 
% 
$$ 
H \vDash \widetilde{f_1} = \widetilde{f_2} 
\;\; \Leftrightarrow \;\; 
H \vDash \widetilde{\varphi}(f_1) = \widetilde{\varphi} (f_2) 
$$ 
and the mapping $\varphi_U$ is a well-defined permutation of $H$ 
(by $\widetilde{\varphi}(a)$ we denote 
the equivalence class of~$\varphi(a)$).

Because 
% 
%-2014dec24-% $$ 
%-2014dec24-% \{i \mid P(\overline f(i))\} \in U \Leftrightarrow 
%-2014dec24-% \{i \mid P(\varphi(\overline f)(i))\} \in U 
%-2014dec24-% $$ 
% 
$$ 
\left\{ 
i \; \left| \; 
P\left(\overline f(i)\right)\right. 
\right\} 
\in U 
\; \Leftrightarrow \; 
\left\{ 
i \; \left| \; 
P\left(\varphi \left( \overline f \right)(i)\right)\right. 
\right\} 
\in U 
$$ 
for every 
% 
%-2014dec24-% $P \in \Sigma$ 
%-2015jan23-% $P \in S$ 
$P \in S,$ 
the mapping 
% 
%---------------------------------------------------------------
%-2015jan01-% $$\varphi_U \colon \prod_U M \to \prod_U M$$ 
$\varphi_U$ 
%---------------------------------------------------------------
% 
preserves all relations from 
% 
%-2014dec24-% $\Sigma$ 
$S$ 
and 
it %-2015jan01-% 
is an automorphism of the structure 
% 
%-2015mar10-% $\prod_U M.$ 
$H.$ 
% 
%-2014dec24-% $\prod_U M$ 
%-2014dec24-% with the signature $\Sigma$. 
% 
Because $V \in U$ the mapping 
% 
%-2015jan01-% $$\varphi \colon \prod_U M \to \prod_U M$$ 
%---------------------------------------------------------------
%-2015jan01-% $$\varphi_U \colon \prod_U M \to \prod_U M$$ 
$\varphi_U$ 
%---------------------------------------------------------------
% 
does not preserve the relation $R$. 
%---------------------------------------------------------------
%-2015jan01-% 
%-2015mar10-% Because $\prod_U M$ is elementary equivalent to $M,$ 
Because $H$ is elementary equivalent to $M,$ 
this is impossible. 
This 
%-2015jan01-%countrudiction 
contradiction 
proves the statement. 
%-2015jan01-% 
%---------------------------------------------------------------
\end{proof} 

%-2015jan23-% We are going to prove: 
We are going to prove now our version of the Svenonius theorem: 

%-2015jan23-% \begin{theorem}(CH) 
%-2015jan23-% Let $S$ be a countable definability space on a universe $A$. 
%-2015jan23-% 
%-2015jan23-% If every permutation on 
%-2015jan23-% % 
%-2015jan23-% %-2014dec24-% $\mathcal F_A$ 
%-2015jan23-% $A^{\mathbb N}$ 
%-2015jan23-% % 
%-2015jan23-% which almost preserves all relations from $S$ 
%-2015jan23-% almost preserves relation $R$ then $R \in S$. 
%-2015jan23-% \end{theorem} 
% 
\begin{theorem}(CH) 
Let $S$ be a countable definability space on a universe~$A$. 

If every permutation on $A^{\mathbb N}$ 
which almost preserves all relations from~$S$ 
almost preserves the relation $R,$ then $R \in S$. 
\end{theorem} 

By $Q^\sigma$ where $Q$ is a formula 
%-2015jan23-% and $\sigma = \pm 1$ 
we denote the formula $Q$ if $\sigma = 1$ 
%-2015jan23-% and $\lnot Q$ if $\sigma = -1$. 
and $\lnot Q$ if $\sigma = 0$. 

%-2015mar10-% %---------------------------------------------------------------
%-2015mar10-% %-2015jan01-% 
%-2015mar10-% We denote by $\Sigma$ a set of names for all elements from $S$ 
%-2015mar10-% and by $M$ the structure $\langle A, \Sigma \cup \{R\}\rangle.$ 
%-2015mar10-% %-2015jan01-% 
%-2015mar10-% %---------------------------------------------------------------
% 
We denote by $\Sigma$ a set of names for all elements from $S,$ 
by $v$ -- the corresponding mapping, 
and by $M$ -- the structure 
$\langle A, \Sigma \cup \{R\}, v \rangle$, and fix a numeration of all elements of $\Sigma$: 
$P_1(\overline x,y),$ 
$P_2(\overline x,y),$ 
$\dots\,.$

The following lemma 
% 
%-2015jan01-% is an analog of 
corresponds to the statement on
$\omega_1$-saturation of ultrapowers 
% 
%-2014dec24-% (e.g. \cite{ult}). 
(\mbox{e.\,g.} \cite{ult}). 
At the same time it corresponds to the inductive step 
in our proof of the Svenonius theorem. 

%-2015jan23-% \begin{lemma} 
%-2015jan23-% If a countable partial mapping 
%-2015jan23-% % 
%-2015jan23-% $$\varphi \colon \mathcal{F} \to \mathcal{F}$$ 
%-2015jan23-% % 
%-2015jan23-% almost preserves a countable definability space~$S,$ 
%-2015jan23-% then for every 
%-2015jan23-% % 
%-2015jan23-% %-2014dec24-% 
%-2015jan23-% $g \not \in Dom(\varphi)$ 
%-2015jan23-% %\mbox{$g \mathop{\not\in} Dom(\varphi)$} 
%-2015jan23-% % 
%-2015jan23-% there exists such 
%-2015jan23-% % 
%-2015jan23-% %-2014dec24-% $h \in \mathcal{F}$ 
%-2015jan23-% \mbox{$h \in \mathcal{F}$} 
%-2015jan23-% % 
%-2015jan23-% that 
%-2015jan23-% % 
%-2015jan23-% %-2014dec24-% $\varphi \cup <g,h>$ 
%-2015jan23-% %-2014dec24-% $\varphi \cup \langle g,h \rangle$ 
%-2015jan23-% \mbox{$\varphi \cup \langle g,h \rangle$} 
%-2015jan23-% % 
%-2015jan23-% almost preserves 
%-2015jan23-% %%BEGIN 
%-2015jan23-% %any relation,  definable in the structure $M$. 
%-2015jan23-% the definability space $S$. 
%-2015jan23-% %END 
%-2015jan23-% \end{lemma} 
% 
\begin{lemma} 
Let $S$ be a countable definability space 
on a universe $A,$ $\mathcal{F} = A^\mathbb{N}.$ 
If a countable partial mapping 
$$\varphi \colon \mathcal{F} \to \mathcal{F}$$ 
almost preserves~$S,$ then for every 
$g \not \in Dom(\varphi)$ there exists such 
\mbox{$h \in \mathcal{F}$} that 
\mbox{$\varphi \cup \langle g,h \rangle$} 
almost preserves the definability space $S.$ 
\end{lemma} 

\begin{proof} 
%-2015jan01-% We use $\Sigma$ as a set of names for $S$ and enumerate 
%-2015jan01-% all formulas in $\Sigma$: 
%-2015jan01-% % 
%-2015jan01-% $$ 
%-2015jan01-% P_1(x_1,\dots,x_{l_1},y), \dots, 
%-2015jan01-% P_m(x_1,\dots,x_{l_m},y), \dots \;. 
%-2015jan01-% $$ 
%-2015jan01-% % 
% 
%-2014dec24-% all 
%-2014dec24-% %%BEGIN 
%-2014dec24-% %elements of 
%-2014dec24-% formulas in 
%-2014dec24-% %END 
%-2014dec24-% $\Sigma$. 
% 
Let us enumerate all elements of  (countable)  $Dom(\varphi)$:
$f_1, \dots, f_n, \dots$ and fix the numeration.
For every $k \in \mathbb{N}$ and 
% 
%-2014dec24-% $$ \overline f = <f_1,\dots,f_l>,$$ 
%-2015jan01-% $$ \overline f = \langle f_1, \dots, f_l \rangle, $$ 
$\overline f = \langle f_1, \dots, f_l \rangle,$ 
$\varphi (\overline f)(k)$ 
denotes 
% 
%-2014dec24-% $$<\varphi (f_1)(k),\dots,\varphi(f_l)(k)>.$$ 
%-2015jan01-% $$ \langle \varphi(f_1)(k), \dots, \varphi(f_l)(k) \rangle.$$ 
$\langle \varphi(f_1)(k), \dots, \varphi(f_l)(k) \rangle.$ 
% 
%-2015jan01-% Denote by 
%-2015jan01-% % 
%-2015jan01-% %-2014dec24-% $M=<\!\!A, \Sigma \cup \{R\}\!\!>$ 
%-2015jan01-% %-2014dec24-% $$M=<\!\!A, \Sigma \cup \{R\}\!\!>$$ 
%-2015jan01-% $$ M= \langle A, \Sigma \cup \{R\} \rangle.$$ 
%-2015jan01-% % 
%-2015jan01-% %-2014dec24-% the countable structure where $\Sigma$ 
%-2015jan01-% %-2014dec24-% is the list of names for relations from $S$.

Let $g \not \in Dom(\varphi)$ be given. 
% 
%-2015jan01-% For a formula 
%-2015jan01-% % 
%-2015jan01-% $$P_i(x_1,\dots,x_{l_i},y)$$ 
%-2015jan01-% % 
%-2015jan01-% and 
%-2015jan01-% %-2014dec24-% a 
%-2015jan01-% $k \in \mathbb{N}$ we denote by 
%-2015jan01-% % 
%-2015jan01-% $$P_i(\overline f(k),g(k))$$ 
%-2015jan01-% % 
%-2015jan01-% the expression 
%-2015jan01-% % 
%-2015jan01-% $$P_i(f_1(k),\dots,f_{l_i}(k),g(k)).$$ 
%-2015jan01-% % 
% 
%-2015jan23-% For any 
For each 
\mbox{$P(x_1,\dots,x_l,y) \in \Sigma$} 
and 
\mbox{$k \in \mathbb{N}$} 
we denote by 
\mbox{$P(\overline f(k),g(k))$} 
the expression 
\mbox{$P(f_1(k),\dots,f_l(k),g(k))$} (we use  $y$ instead of, for example  $x_{l+1}$, for the convenience of our notations).

\sloppy %-2015jan23-% 

%---------------------------------------------------------------
%-2015jan01-% 
%We enumerate all elements of $\Sigma$: 
%$P_1(\overline x,y),$ 
%$P_2(\overline x,y),$ 
%$\dots\,.$ 
%-2015jan01-% 
%---------------------------------------------------------------
% 
For each 
\mbox{$k \in \mathbb{N}$} 
let us take 
\mbox{$\sigma_{i,k}=1$} 
if 
% 
%-2015jan01-% $$M \vDash P_i(\overline f(k),g(k))$$ 
\mbox{$M \vDash P_i(\overline f(k),g(k))$} 
% 
%-2014dec24-% and $\sigma_i=-1$ if 
and 
%-2015jan23-% \mbox{$\sigma_i \mathop{=} -1$} 
\mbox{$\sigma_{i,k} \mathop{=} 0$} 
if 
% 
%-2015jan01-% $$M \vDash \lnot P_i(\overline f(k),g(k)).$$ 
\mbox{$M \vDash \lnot P_i(\overline f(k),g(k)).$} 
% 
%-2014dec24-% and define $m(k)$ as 
% 
% 
% 
We define $m(k)$ as maximal element of the set
% 
%-2014dec24-% $$ 
%-2014dec24-% \max \{m \leqslant k | 
%-2014dec24-% M \vDash (\exists y) 
%-2014dec24-% (\bigwedge_{i=1}^m P_i^{\sigma_i} 
%-2014dec24-% (\varphi (\overline f)(k),y))\} 
%-2014dec24-% $$ 
$$ 
%\mathrel{\max} 
\left\{ 
m \leqslant k \; \; \left| \; \; 
M \vDash (\exists y) 
\left( 
\bigwedge_{i=1}^m P_i^{\sigma_{i,k}} 
\biggl( % \left( 
\varphi \left( \overline f \right)(k),y 
\biggr) % \right) 
\right) 
\right. 
\right\} 
$$ 
and 
$m(k)=0$ if the set 
% 
%-2014dec24-% $$ 
%-2014dec24-% \{m \leqslant k | M \vDash 
%-2014dec24-% (\exists y)(\bigwedge_{i=1}^m 
%-2014dec24-% P_i^{\sigma_i}(\varphi (\overline f)(k),y))\} 
%-2014dec24-% $$ 
%$$ 
%%\mathrel{\max} 
%\left\{ m \leqslant k \; \; \left| \; \; M \vDash (\exists y) \left( \bigwedge_{i=1}^m P_i^{\sigma_{i,k}} \biggl( % \left( \varphi \left( \overline f \right)(k),y \biggr) % \right) \right) \right. \right\} 
%$$ 
% 
is empty. 
% 
% 
%-2014dec24-% Define 
We define 
a function 
% 
%-2015jan01-% $$h \colon \mathbb{N} \to A$$ 
\mbox{$h \colon \mathbb{N} \to A$} 
% 
%-2015jan01-% such that 
so that 
% 
%-2014dec24-% $$ 
%-2014dec24-% M \vDash \bigwedge_{i=1}^{m(k)} 
%-2014dec24-% P_i^{\sigma_i} 
%-2014dec24-% (\varphi (\overline f)(k),h(k)) 
%-2014dec24-% $$ 
$$ 
M \vDash \bigwedge_{i=1}^{m(k)} 
P_i^{\sigma_{i,k}} 
\biggl( \varphi \left( \overline f \right)(k),h(k) \biggr) 
$$ 
if 
% 
%-2014dec24-% $m(k) > 0$, 
$m(k) \mathop{>} 0$, \, 
$h(k)$ is an arbitrary element of $A$ 
if $m(k)=0$.

We need to show that 
% 
%-2014dec24-% $$ 
%-2014dec24-% \{i | 
%-2014dec24-% P_j(f_1(i),\dots,f_{l_j}(i),g(i)) \not \equiv 
%-2014dec24-% P_j(\varphi(f_1)(i),\dots,\varphi(f_{l_j})(i),h(i) 
%-2014dec24-% )\} 
%-2014dec24-% $$ 
% 
%-2015jan01-% $$ 
%-2015jan01-% \left\{ 
%-2015jan01-% i \; \; \left| \; \; 
%-2015jan01-% P_j \bigl( f_1(i),\dots,f_{l_j}(i),g(i) \bigr) 
%-2015jan01-% \not \equiv 
%-2015jan01-% P_j \bigl( \varphi(f_1)(i),\dots,\varphi(f_{l_j})(i),h(i) \bigl) 
%-2015jan01-% \right. 
%-2015jan01-% \right\} 
%-2015jan01-% $$ 
% 
$$ 
\left\{ 
i \; \; \left| \; \; 
P_j \bigl(\overline f(i),g(i) \bigr) 
\not \equiv 
P_j \bigl( \varphi(\overline f)(i),h(i) \bigl) 
\right. 
\right\} 
$$ 
% 
%-2015jan01-% is finite for each $j$. 
is finite for each $P_j \in \Sigma.$ 
Consider all formulas 
$$ Q_{\overline \tau}(\overline x,y) = 
\bigwedge_{i=1}^j P_i^{\tau_i}(\overline x,y) 
$$ 
% 
%-2015jan23-% for all collections 
for all tuples 
% 
%-2014dec24-% $\overline\sigma = 
%-2014dec24-% (\sigma_1,\dots\sigma_j) \in \{0,1\}^\mathbb{N}$. 
% 
%-2015jan01-% $$ \overline\sigma = 
%-2015jan01-% (\sigma_1,\dots\sigma_j) \in \{0,1\}^\mathbb{N}.$$ 
% 
%-2015jan23-% $\overline\sigma = 
%-2015jan23-% (\sigma_1,\dots\sigma_j) \in \{0,1\}^\mathbb{N}.$ 
$\overline\tau = 
(\tau_1,\dots\tau_j) \in \{0,1\}^j.$ 
The mapping $\varphi$ almost preserves all relations 
% 
%-2014dec24-% $ (\exists y)Q_{\overline \sigma}(\overline x,y) $ 
% 
%-2015jan01-% $$ (\exists y)Q_{\overline \sigma}(\overline x,y) $$ 
$(\exists y)Q_{\overline \tau}(\overline x,y).$ 
% 
%-2015jan01-% so there is such $n_0$ that 
%-2015jan23-% So there is such $n_0$ that 
So, there is such $n_0$ that 
% 
%-2014dec24-% $ (\exists y) Q_{\overline\sigma} (\overline f(i),y) 
%-2014dec24-% \equiv 
%-2014dec24-% (\exists y)Q_{\overline\sigma} (\varphi(\overline f)(i),y)$ 
$$ 
\mathop{(\exists y)} Q_{\overline\tau} 
\left( \overline f(i),y \right) 
\equiv 
\mathop{(\exists y)} Q_{\overline\tau} 
\left( \varphi(\overline f)(i),y \right) 
$$ 
for every tuple $\overline\tau$ and 
% 
%-2014dec24-% $i > n_0$. 
$i \mathop{>} n_0$. 
If 
% 
%-2014dec24-% $k>\max\{j,n_0\}$ 
$k \mathop{>} \max \, \{ j, n_0 \}$ 
then by definition of $m(k)$ we have 
% 
%-2014dec24-% $m(k) > j$ 
$m(k) \mathop{>} j$ 
so 
%-2014dec24-% $$ 
%-2014dec24-% P_j(f_1(k),\dots,f_{l_j}(k),g(k)) \equiv 
%-2014dec24-% P_j(\varphi(f_1)(k),\dots,\varphi(f_{l_j})(k),h(k)). 
%-2014dec24-% $$ 
% 
%-2015jan01-% $$ 
%-2015jan01-% P_j \bigl( f_1(k), \dots, f_{l_j}(k), g(k) \bigr) 
%-2015jan01-% \equiv 
%-2015jan01-% P_j \bigl( \varphi(f_1)(k), \dots, 
%-2015jan01-% \varphi(f_{l_j})(k), h(k) \bigr). 
%-2015jan01-% $$ 
% 
%-2015jan01-% $$ 
%-2015jan01-% P_j \bigl( \overline f(k),g(k) \bigr) 
%-2015jan01-% \equiv 
%-2015jan01-% P_j \bigl( \varphi(\overline f)(k),h(k) \bigr). 
%-2015jan01-% $$ 
% 
\hfill 
$ 
P_j \bigl( \overline f(k),g(k) \bigr) 
\equiv 
P_j \bigl( \varphi(\overline f)(k),h(k) \bigr). 
$ 
\hfill 
\end{proof} 

\bigskip %-2015jan01-% 
%-2014dec24-% \textbf{Proof of the theorem.} 
\textbf{Proof of the theorem} 

%-2014dec24-% \begin{proof} 
%-2015jan01-% Denote by 
%-2015jan01-% %-2014dec24-% $M=<\!\!A, \Sigma \cup \{R\}\!\!>$ 
%-2015jan01-% $ M {=} \langle A, \Sigma \cup \{R\} \rangle$ 
%-2015jan01-% the countable structure where $\Sigma$ is the list of 
%-2015jan01-% names for relations from~$S$. 
% 
%-2015jan01-% Assume that $R \not \in S$ and construct 
We assume that $R \not \in S$ and will construct 
a permutation $\varphi$ of $\mathcal{F}$ 
that almost preserves all relations in $S$ but does not 
almost preserve $R$. 

%-2015mar10-% If $f,g \in \mathcal{F}$ then $f \approx g$ means that 
For $f,g \in \mathcal{F}$ we write $f \approx g$ when 
% 
%-2015jan01-% $$ 
%-2015jan01-% \{i \mid f(i) \ne g(i)\} 
%-2015jan01-% $$ 
% 
\mbox{$\{i \mid f(i) \ne g(i)\}$} 
is finite, 
by $[f]$ we denote the corresponding equivalence class of $f$. 
% 
%-2015jan01-% According CH 
Using CH 
we can order $\mathcal{F}$ as~$\omega_1$. 
For each ordinal 
$\alpha < \omega_1$ 
we will construct such a partial bijection 
% 
%-2015jan01-% $$\varphi_\alpha \colon \mathcal{F} \to \mathcal{F}$$ 
\mbox{$\varphi_\alpha \colon \mathcal{F} \to \mathcal{F}$} 
that 

(i) 
$Dom(\varphi_\alpha)$ 
is countable 

(ii) 
$\varphi_\alpha \subseteq \varphi_\beta$ 
if 
$\alpha < \beta$ 

(iii) 
if 
$f \in Dom(\varphi_\alpha)$ 

\qquad 
then $[f] \subseteq Dom(\varphi_\alpha)$ 
and 
$\varphi_\alpha([f])=[\varphi_\alpha(f)]$ 

%-2015jan01-% (iv) 
%-2015jan01-% $\varphi_\alpha$ 
%-2015jan01-% almost preserves every relation definable in 
%-2015jan01-% % 
%-2015jan01-% %-2014dec24-% $<\!\!A, \Sigma\!\!>$. 
%-2015jan01-% 
%-2015jan01-% \qquad $\langle A, \Sigma\rangle$. 
% 
(iv) 
$\varphi_\alpha$ 
almost preserves every relation 
%-2015jan01-% definable in $\langle A, \Sigma\rangle$. 
from $S.$                                     %-2015jan01-% 

\smallskip 
\noindent 
%-2015jan01-% We set 
Then we shall set 
$\varphi = \bigcup_{\alpha < \omega _1} \varphi_\alpha$. 
 
\medskip 
$\boldsymbol{\alpha=0}$. 
%-2015jan01-% We enumerate 
%-2015jan01-% all relations from $S$: 
%-2015jan01-% $$ 
%-2015jan01-% P_1(\overline x), \dots, P_m(\overline x), \dots 
%-2015jan01-% $$ 
%-2015jan01-% % 
%-2015jan01-% %-2014dec24-% all relations from $S$ 
%-2015jan01-% % 
%-2015jan01-% and note that 
Let us note that 
$$ 
M \vDash 
(\exists \overline a)(\exists \overline b) 
\left( 
\Bigl( R(\overline a) \not \equiv R(\overline b) \Bigr) 
\land 
\bigwedge_{i=1}^m 
\biggl( 
P_i(\overline a) 
\equiv 
P_i(\overline b) 
\biggr) 
\right) 
$$ 
for all 
$m \in \mathbb{N}$. 
Otherwise 
$$ 
M \vDash 
(\forall \overline a) 
(\forall \overline b) 
\left( 
\bigwedge_{i=1}^m 
\left( P_i(\overline a) \equiv P_i(\overline b) \right) 
\to 
\left( R(\overline a) \equiv R(\overline b) \right) 
\right) 
$$ 
% 
%-2015jan23-% for some $m \in \mathbb{N}$ so 
%-2015jan23-% $R(\overline x)$ 
%-2015jan23-% is equivalent 
for some $m \in \mathbb{N},$ so, 
$R(\overline x)$ 
would be equivalent 
to a disjunction of finite number of 
% 
%-2015jan01-% $$ 
%-2015jan01-% \bigwedge_{i=1}^m P_i^{\sigma_i}(\overline x) 
%-2015jan01-% $$ 
%-2015jan01-% formulas. 
formulas 
\mbox{$\bigwedge_{i=1}^m P_i^{\sigma_i}(\overline x).$} 

\medskip %-2015jan01-% 
For each $m \in \mathbb{N}$ let 
$\overline a(m)$ and $\overline b(m)$ 
be tuples such that 
% 
%-2014dec24-% $$ 
%-2014dec24-% M \vDash 
%-2014dec24-% ( 
%-2014dec24-% (R(\overline a(m)) \not \equiv  R(\overline b(m))) 
%-2014dec24-% \land 
%-2014dec24-% \bigwedge_{i=1}^m 
%-2014dec24-% (P_i(\overline a(m)) \equiv P_i(\overline b(m))) 
%-2014dec24-% ). 
%-2014dec24-% $$ 
\begin{equation} 
M \vDash 
\left( 
\Bigl( 
R(\overline a(m)) \not \equiv  R \left( \overline b(m) \right) 
\Bigr) 
\land 
\bigwedge_{i=1}^m 
\Bigl( 
P_i(\overline a(m)) \equiv P_i \left( \overline b(m) \right) 
\Bigr) 
\right). 
\end{equation}
We define functions 
% 
%-2015jan01-% $$ g_1,\dots,g_n,h_1,\dots,h_n \in \mathcal{F} $$ 
% 
$g_1,$ $\dots,$ $g_n,$ $h_1,$ $\dots,$ 
\mbox{$h_n \in \mathcal{F}$} 
% 
%-2015jan23-% where $n$ is arity of $R$ as 
where $n$ is the arity of~$R$ as 
$$ 
g_t(i)=(\overline a(i))_t, 
\quad 
h_t(i)=(\overline b(i))_t, 
%-2015jan23-% \quad \varphi_0(g_i)=h_i 
$$ 
% 
%-2015jan01-% and expand $\varphi_0$ so that 
and define $\varphi_0$ so that 
% 
%$\varphi_0(g_i)=h_i$ \, %-2015jan23-% 
% 
%-2015jan01-% $$\varphi_0([g_i])=[h_i].$$ 
\mbox{$\varphi_0([g_i])=[h_i]$} is an arbitrary bijection $[g_i] \to [h_i]$
\, %-2015jan23-% 
$i=1,$ $\dots,$ $n.$ 

Let us remind, that we consider structures with equality, and $(x_i=x_j)$ is $P_k(\bar x)$ for some $k$. So $g_m \approx g_l \Leftrightarrow h_m \approx h_l$ for all $m,l$ and the definition of $\varphi_0$ is correct.

%-----------------------------
%-2015jan01-% 
So, the condition (iii) holds. 
%-2015jan01-% 
%-----------------------------
% 
Because $[f]$ is countable 
%-2015jan23-% for any $f$  %-2015jan01-% 
for every $f,$  %-2015jan01-% %-2015jan23-% 
the condition (i) holds. 
%-2015jan01-% %-2014dec24-% It's 
%-2015jan01-% It is 
%-2015jan01-% easy to note that 
It is easy to see due to definition (1) that 
% 
%-2014dec24-% $$ 
%-2014dec24-% M \vDash P_i(g_1(j), \dots, g_n(j)) 
%-2014dec24-% \equiv 
%-2014dec24-% P_i(\varphi_0(g_1)(j), \dots, \varphi_0(g_n)(j)) 
%-2014dec24-% $$ 
$$ 
M \vDash 
P_i \bigl( g_1(j), \dots, g_n(j) \bigr) 
\equiv 
P_i \bigl(\varphi_0(g_1)(j), \dots, \varphi_0(g_n)(j) \bigr) 
$$ 
% 
%-2015jan01-% for all $j > i$ and 
%-2015jan23-% for all $j > i$ 
for all $j > i,$ 
so the condition (iv) holds. 
We see also, that 
% 
%-2014dec24-% $$ 
%-2014dec24-% M \vDash 
%-2014dec24-% R(g_1(j),\dots,g_n(j)) \not \equiv 
%-2014dec24-% R(\varphi_0(g_1)(j),\dots,\varphi_0 (g_n)(j)) 
%-2014dec24-% $$ 
\begin{equation}  
M \vDash 
R \bigl( g_1(j), \dots, g_n(j) \bigr) 
\not \equiv 
R \bigl( \varphi_0(g_1)(j), \dots, \varphi_0 (g_n)(j) \bigr) 
\end{equation} 
for all $j$. 
%---------------------------------------------------------------
%-2015jan01-% 
We shall have this condition for all 
$\varphi_\alpha,$ $\alpha>0$ . 
%-2015jan01-% 

We shell prove (ii) for all $\alpha < \beta$
 and conclude that $\varphi_0 \subset \varphi_\alpha$ for $\alpha > 0$. So, we shall have the statement (2) for all $\varphi_\alpha$ instead of $\varphi_0$. 
%---------------------------------------------------------------

\medskip 
$\boldsymbol{0<\alpha<\omega _1}$. 
% 
%-2015jan01-% Denote by 
%-2015jan01-% % 
%-2015jan01-% $$\varphi'=\bigcup_{\beta < \alpha} \varphi_\beta$$ 
%-2015jan01-% % 
% 
Let us denote by 
$\varphi'$ 
the countable mapping 
$\bigcup_{\beta < \alpha} \varphi_\beta$ 
and by $g$ the 
% 
%-2014dec24-% item 
element 
of $\mathcal{F}$ with the index $\alpha$. 
%---------------------------------------------------------------
%-2015jan01-% 
%-2015jan23-% We define 
We shall define 
$\varphi_\alpha \supset \varphi'$. 
%-2015jan01-% 
%---------------------------------------------------------------

\sloppy %-2015jan23-% 

%-2015jan23-% Suppose that $g \not \in Dom(\varphi')$. 
Let us suppose that $g \not \in Dom(\varphi')$. 
% 
%-2015jan01-% We define 
%-2015jan01-% $\varphi_\alpha \supset \varphi'$. 
% 
%-2015jan01-% %%ADDED 
%-2015jan01-% The set 
%-2015jan01-% % 
%-2015jan01-% %-2014dec24-% $$\{\beta | \beta < \alpha \}$$ 
%-2015jan01-% $$\{\beta \mid \beta < \alpha \}$$ 
%-2015jan01-% % 
%-2015jan01-% is countable, so we can 
%-2015jan01-% %ADDED 
% 
We %-2015jan01-% 
use lemma 1 to find 
% 
%-2015jan01-% corresponding $h \in \mathcal{F}$. 
$h \in \mathcal{F}$ corresponding to $g.$ 
% 
%-2014dec24-% According condition (iii) 
According to condition (iii) 
% 
%-------- 2015mar10 --------% 
for $\beta < \alpha$: 
%-------- 2015mar10 --------% 
% 
%-2015jan01-% $$[g] \cap Dom(\varphi') = \varnothing$$ 
$$[g] \cap Dom(\varphi') = \varnothing.$$ 
% 
%-2014dec24-% so it's easy to see that 
%-2015jan01-% so it is easy to see that 
% 
%---------------------------------------------------------------
%-2015jan01-% 
At the same time 
$\varphi' \cup \{\langle g,h \rangle\}$ 
almost preserves $=,$ so 
%-2015jan01-% 
%---------------------------------------------------------------
% 
%-2015jan01-% $$[h] \cap Im(\varphi') = \varnothing$$ 
$$[h] \cap Im(\varphi') = \varnothing,$$ 
% 
%-2015jan01-% as well 
% 
and we can define 
$\varphi_\alpha(g)=h$. 
Similarly we can define 
$\varphi_\alpha^{-1}(g)$ 
if 
%-2015jan23-% $g \not \in Im(\varphi')$. 
\mbox{$g \not \in Im(\varphi')$.} 
% 
%-2015jan01-% We extend 
We extend now 
$\varphi_\alpha$ so that 
$\varphi_\alpha([g]) = [\varphi_\alpha(g)]$ 
and 
%-2015jan23-% $\varphi_\alpha^{-1}([g]) = [\varphi_\alpha^{-1}(g)]$. 
\mbox{$\varphi_\alpha^{-1}([g]) = [\varphi_\alpha^{-1}(g)],$} 
and we keep the condition (iii)  %-2015jan23-% 
valid for $\varphi_\alpha.$      %-2015jan23-% 
% 
%---------------- 2015mar10 ----------------% 
Other conditions are evident. 
%---------------- 2015mar10 ----------------% 
% 
\hfill $\Box$ 
%-2014dec24-% \end{proof} 

%--------------------------------------------------------------------------
%---------------------------- 2015jan23 ----------------------------%
\medskip 
Our proof of the theorem uses induction essentially. 
So, we do not know how to prove it without CH.  
%---------------------------- 2015jan23 ----------------------------%
%--------------------------------------------------------------------------

%==========================================================================
%//////////////////////////////////////////////////////////////////////////
%==========================================================================

%--------------------------------------------------------------------------
%---------------------------- 2015jan25 ----------------------------%
\bigskip 
{\bf 5. Concluding remarks and open problems} 
\nopagebreak 
\bigskip 
%---------------------------- 2015jan25 ----------------------------%
%--------------------------------------------------------------------------

If $G$ is a group of permutations of a set $A$, 
then we say that $G$ is \emph{closed} if for every $g \in Sym(A)$, if 
%the following holds: suppose that $g \in Sym(A)$ and 
for every tuple $\overline a \in A$ there is 
$h \in G$ such that 
% 
%-2015jan01-% $g(\overline a)=h(\overline a)$ 
$g(\overline a)=h(\overline a),$ 
then $g \in G$. 

%-2014dec24-% It's well known (e.g. \cite{mac}) 
It is well known (\mbox{e.\,g.} \cite{mac}) 
that for every definability space $S$ 
the corresponding group $G_S$ is closed. 
% 
%-2015jan01-% From other hand 
On the other hand 
each closed subgroup of $Sym(A)$ 
% 
%-2015jan01-% is the group of permutations, 
is the group of all permutations, 
% 
%-2015jan01-% preserving some definability space: 
%-2015jan01-% with a subgroup $G$ of $Sym(A)$ 
which preserve some definability space. 
% 
%-2015jan25-% With any subgroup $G$ of $Sym(A)$ 
With every subgroup $G$ of $Sym(A)$ 
we can associate 
% 
%-2015mar10-% a \emph{canonical basis}, 
a countable \emph{canonical basis}, 
% 
%-2014dec24-% i.e. 
\mbox{i.\,e.} 
for each $n < \omega$ 
and each orbit of $G$ in $A^n$ 
%-2015jan25-% we choose an $n-$ary relation. 
we choose an $n$-ary relation. 
So the next statement holds: 

\begin{statement} 
\cite{mac}                       %-2015jan01-% 
%-2015jan01-% Subgroup $G$ 
A subgroup $G$ 
of the group $Sym(A)$ 
is the group of automorphisms 
of some countable definability space 
iff $G$ is closed. 
\end{statement} 

In our theorem we consider the group $Sym(\mathcal{F})$ 
and with each definability space $S$ 
we associate the subgroup 
$G^*_S \subseteq Sym(\mathcal{F})$ 
of permutations almost preserving relations from $S$. 

\begin{que} $  $ 

%-2015jan01-% (i) Could we 
(i) Can we 
describe all subgroups of 
$Sym(\mathcal{F})$ 
%-2014dec24-% of form $G^*_S$? 
of the form $G^*_S$? 

%-2015jan01-% (ii) Could we describe all subgroups of 
%-2015jan01-% $Sym(\mathcal{F})$ of form $G^*_S$ 
%-2015jan01-% for a definability space $S$ with a finite basis? 
%-2015jan01-% for decidable structures 
%-2015jan01-% % 
%-2015jan01-% %-2014dec24-% $<\!\!A,\Sigma\!\!>$? 
%-2015jan01-% $ \langle A,\Sigma \rangle$? 
% 
(ii) Can we describe all subgroups of 
$Sym(\mathcal{F})$ of form $G^*_S$ 
for definability spaces $S$ with finite basis? 
for decidable structures 
$\langle A,\Sigma \rangle$? 
\end{que} 

Our theorem states that $R$ is not definable in $S$ 
if there is no permutation $\varphi \in Sym(\mathcal{F})$ 
that almost preserves $S$ but not $R.$ 
But the group $Sym(\mathcal{F})$ is too big. 
Can we limit our search of $\varphi$ with a "natural" 
smaller subgroup of $Sym(\mathcal{F})$? 

%-2014dec24-% —ÙÓÏÛÎËÓ‚‡ÌÌ‡ˇ ÚÂÓÂÏ‡ ÛÚ‚ÂÊ‰‡ÂÚ, 
%-2014dec24-% ˜ÚÓ ÓÚÌÓ¯ÂÌËÂ $R$ ÓÔÂ‰ÂÎËÏÓ 
%-2014dec24-% ˜ÂÂÁ ÓÚÌÓ¯ÂÌËˇ ÒË„Ì‡ÚÛ˚ $\Sigma$ 
%-2014dec24-% ÂÒÎË ÌÂ ÒÛ˘ÂÒÚ‚ÛÂÚ ÍÓÌÚÔËÏÂ‡: ÔÂÂÒÚ‡ÌÓ‚ÍË 
%-2014dec24-% $\varphi \in Sym(\mathcal{F})$, 
%-2014dec24-% ÔÓ˜ÚË ÒÓı‡Ìˇ˛˘ÂÈ ‚ÒÂ ÓÚÌÓ¯ÂÌËˇ ËÁ $\Sigma$ 
%-2014dec24-% Ë ÌÂ ÔÓ˜ÚË ÒÓı‡Ìˇ˛˘ÂÈ ÓÚÌÓ¯ÂÌËÂ $R$. 
%-2014dec24-% √ÛÔÔ‡ $Sym(\mathcal{F})$ 
%-2014dec24-% ˜ÂÁ‚˚˜‡ÈÌÓ ·ÓÎ¸¯‡ˇ, 
%-2014dec24-% ÏÓÊÌÓ ÎË ÔË ÔÓËÒÍÂ ÍÓÌÚÔËÏÂ‡ 
%-2014dec24-% Ó„‡ÌË˜ËÚÒˇ ÌÂ·ÓÎ¸¯ÓÈ ÂÒÚÂÒÚ‚ÂÌÌÓÈ 
%-2014dec24-% ÔÓ‰„ÛÔÔÓÈ „ÛÔÔ˚ $Sym(\mathcal{F})$? 

\begin{que} 
%-2014dec24-% —Û˘ÂÒÚ‚ÛÂÚ ÎË Ú‡Í‡ˇ ÂÒÚÂÒÚ‚ÂÌÌ‡ˇ ÔÓ‰„ÛÔÔ‡ 
Is there any "natural" subgroup 
% 
%-2014dec24-% $H \subset Sym(\mathcal{F})$, 
%-2014dec24-% \mbox{$H \subset Sym(\mathcal{F}),$} 
$$H \subset Sym(\mathcal{F}),$$ 
that 
for any $S_1,$ $S_2$ 
% 
%-2014dec24-% $H \cap G^*_{S_1} = H \cap G^*_{S_2} 
%-2014dec24-% \Leftrightarrow S_1 = S_2$ 
% 
$$ 
H \cap G^*_{S_1} = H \cap G^*_{S_2} 
\; \; \Leftrightarrow \; \; 
S_1 = S_2 \; \; ? 
$$ 
\end{que} 

%-2015jan01-% From simple cardinality arguments 
%-2015jan01-% follows that there is a subgroup 
From simple cardinality arguments 
we conclude that there is a subgroup 
% 
%-2014dec24-% $H \subset Sym(\mathcal{F}), |H|=2^{\aleph_0}$ 
$$H \subset Sym(\mathcal{F}), \quad |H|=2^{\aleph_0}$$ 
with the required property, 
but it is desirable to get an explicit description. 

\textbf{Remark.} 
With a structure 
% 
%-2014dec24-% $M=<\!\!A, \Sigma\!\!>$ 
$ M = \langle A,\Sigma \rangle$ 
we can associate a Boolean-valued structure 
% 
%-2014dec24-% $M_\Sigma=<B,<\mathcal F_A, \Sigma>,\Psi>$ 
% 
% 
%-2014dec24-% $$ 
%-2014dec24-% M_\Sigma = 
%-2014dec24-% \langle 
%-2014dec24-% B, \langle \mathcal F_A, \Sigma \rangle, \Psi 
%-2014dec24-% \rangle 
%-2014dec24-% $$ 
% 
% 
$$ 
M_\Sigma = 
\bigl\langle 
%-2015jan25-% B, \langle \mathcal F_A, \Sigma \rangle, \Psi 
B, \langle A^{\mathbb N}, \Sigma \rangle, \Psi 
\bigr\rangle 
$$ 
such that 

(1) 
% 
%-2015jan01-% $B$ is Boolean algebra 
$B$ is a Boolean algebra 
$2^\mathbb{N} / \approx $, 
where the equivalence $\approx$ is defined as 
% 
%-2014dec24-% $A_1 \approx A_2 \Leftrightarrow 
%-2014dec24-% ( 
%-2014dec24-% (A_1 \setminus A_2) \cup 
%-2014dec24-% (A_2 \setminus A_1)$ is finite). 
% 
$$ 
A_1 \approx A_2 
\; \; \Leftrightarrow \; \; 
\bigl( 
(A_1 \setminus A_2) \cup (A_2 \setminus A_1) 
\; \text{ is finite} 
\bigr) 
. 
$$ 

(2) 
$\Psi$ is the operation which assigns to each formula 
% 
%-2015jan01-% $Q(x_0,\dots,x_n)$ 
$Q(x_1,\dots,x_n)$ 
of $\Sigma$ with free variables among 
% 
%-2015jan01-% $x_0,\dots,x_n$ 
$x_1,\dots,x_n$ 
a function 
% 
%-2015jan01-% $A^{n+1} \to B$, 
$A^n \to B$, 
% 
%-2015jan01-% whose value at 
%-2015jan01-% $f_0,\dots,f_n$ 
%-2015jan01-% is denoted 
%-2015jan01-% $\Psi(Q(f_0,\dots,f_n))$. 
% 
with value at 
$f_1,$ $\dots,$ $f_n$ 
denoted by 
$\Psi(Q(f_1,\dots,f_n))$. 

\smallskip 
(2.1) 
For a relation 
$P \in \Sigma$ 
we set 
% 
%-2014dec24-% $$ 
%-2014dec24-% \Psi(P(f_0, \dots, f_n)) = 
%-2014dec24-% \{ i | P(f_0(i), \dots, f_n(i)) \} / \approx 
%-2014dec24-% ; 
%-2014dec24-% $$ 
% 

\smallskip 
%-2015jan01-% $ 
%-2015jan01-% \Psi(P(f_0, \dots, f_n)) = 
%-2015jan01-% \{ i \mid P(f_0(i), \dots, f_n(i)) \} / \approx 
%-2015jan01-% ; 
%-2015jan01-% $ 
% 
$ 
\Psi(P(f_1, \dots, f_n)) = 
\{ i \mid P(f_1(i), \dots, f_n(i)) \} / \approx 
; 
$ 

\smallskip 
$\Psi(P \lor Q) = \Psi(P) \lor \Psi(Q);$ 

\smallskip 
$\Psi(\lnot P) = -\Psi(P)$ (the complement of $\Psi(P)$); 

\smallskip 
$\Psi(\exists x P(\overline u,x)) = 
\bigvee 
\{ \Psi(P(\overline u,f)) \mid f \in \mathcal{F}\}.$ 

\smallskip 
%-2015jan01-% Note, that the supremum exists in this case. 
Note, that the supremum exists. 

\medskip 
% 
%-2014dec24-% It's easy to see that 
It is easy to see that 
% 
%-2014dec24-% $$ 
%-2014dec24-% \Psi(Q(f_0, \dots, f_n)) = 
%-2014dec24-% \{i | Q(f_0(i), \dots, f_n(i))\} / \approx 
%-2014dec24-% $$ 
% 
%-2015jan01-% $$ 
%-2015jan01-% \Psi(Q(f_0, \dots, f_n)) = 
%-2015jan01-% \{i \; \mid \; Q(f_0(i), \dots, f_n(i))\} / \approx 
%-2015jan01-% $$ 
%-2015jan01-% %  
$$ 
\Psi(Q(f_1, \dots, f_n)) = 
\{i \; \mid \; Q(f_1(i), \dots, f_n(i))\} / \approx 
$$ 
% 
%-2015jan01-% for every formula $Q$ in $\Sigma$, 
%-2015jan01-% and every permutation $\varphi$ of $\mathcal F_A$ 
%-2015jan01-% almost preserving relations from $\Sigma$ 
%-2015jan01-% is the automorphism of $M_\Sigma$. 
% 
%-2015mar10-% is an automorphism of $M_\Sigma$ 
%-2015mar10-% for any formula $Q$ in $\Sigma$, 
%-2015mar10-% %-2015jan25-% and any permutation $\varphi$ of $\mathcal F_A$ 
%-2015mar10-% and any permutation $\varphi$ of $A^{\mathbb N}$ 
%-2015mar10-% almost preserving relations from $\Sigma.$ 
% 
for every formula $Q$ in $\Sigma,$ 
and every permutation $\varphi$ of $A^{\mathbb N}$ 
almost preserving relations from $\Sigma$ 
is an automorphism of $M_\Sigma.$ 
So our theorem may be reformulated 
% 
%-2014dec24-% in terms of Boolean valued structure as 
in terms of Boolean-valued structures as 

\begin{statement} 
(CH). 
Let 
% 
%-2014dec24-% $M=<\!\!A, \Sigma \cup \{R\}\!\!>$ 
$ M = \langle A, \Sigma \cup \{R\} \rangle$ 
be a countable structure. 
The following conditions are equivalent: 

(1) 
A relation $R$ is definable in 
% 
%-2014dec24-% $<\!\!A, \Sigma\!\!>$. 
$ \langle A, \Sigma \rangle$. 

(2) 
Every automorphism of the structure $M_\Sigma$ 
% 
%-2015jan01-% is the automorphism of the structure 
is an automorphism of the structure 
$M_{\Sigma \cup \{R\}}$. 
\end{statement} 

Let us mention the paper \cite{bool} where 
Boolean-valued structures were used to describe 
first-order definability also: 

\begin{statement}\!\!\cite{bool} 

Let $T$ be any first order theory. 
There exists a Boolean valued structure 
$\mathfrak{M}$ such that 

(i) 
$\mathfrak{M}$ is a conservative model of $T$, 
in the sense that 
$\mathfrak{M} \vDash \varPhi$ 
iff 
$T \vdash \varPhi$, 
for each sentence $\varPhi$. 

(ii) 
Any predicate which is invariant under all 
automorphisms of $\mathfrak{M}$ is definable. 
\end{statement} 

%-2014dec24-% œÂ‰‚‡ËÚÂÎ¸Ì‡ˇ ‚ÂÒËˇ ÚÂÍÒÚ‡ 
%-2014dec24-% ·˚Î‡ ‡ÁÏÂ˘ÂÌ‡ ‚ \cite{se1}, 
%-2014dec24-% ÂÁÛÎ¸Ú‡Ú ·˚Î ÒÓÓ·˘ÂÌ ‚ \cite{se2}. 

%==========================================================================
%//////////////////////////////////////////////////////////////////////////
%==========================================================================

%--------------------------------------------------------------------------
%---------------------------- 2015jan25 ----------------------------%
%-2015jan25-% {6. Õ‡ Ò. 16} 

\end{document}